\documentclass[11pt,a4paper]{amsart}
\usepackage[T1]{fontenc}
\usepackage{lmodern}
\usepackage{amsmath,amssymb,amsthm,mathtools}
\usepackage[margin=30mm]{geometry}
\usepackage{enumitem}
\usepackage[dvipdfmx,hidelinks]{hyperref}
\usepackage{url}
\hypersetup{
  pdftitle={The Z/2-equivariant Seiberg--Witten Floer cohomology
  of the (2,1)-cable detects the unknot},
  pdfauthor={Nobuo Iida and Tatsumasa Suzuki}
}
\numberwithin{equation}{section}
\newtheorem{theorem}{Theorem}[section]
\newtheorem{proposition}[theorem]{Proposition}
\newtheorem{lemma}[theorem]{Lemma}
\newtheorem{corollary}[theorem]{Corollary}
\theoremstyle{definition}

\theoremstyle{remark}

\newcommand{\Z}{\mathbb Z}
\newcommand{\Q}{\mathbb Q}

\newcommand{\Spinc}{\operatorname{Spin}^{c}}

\newcommand{\cab}{\operatorname{cab}}

\newcommand{\s}{\mathfrak s}

\title[Equivariant Floer cohomology of a cable]{The $\Z/2$-equivariant Seiberg--Witten Floer cohomology of the $(2, 1)$-cable is an unknot detector}

\author{Nobuo Iida}
\address{Kavli Institute for the Physics and Mathematics of the Universe (WPI), The University of Tokyo, 5-1-5 Kashiwanoha, Kashiwa, Chiba 277-8583, Japan}
\email{iidanobuo1224@g.ecc.u-tokyo.ac.jp}

\author{Tatsumasa Suzuki}
\address{Graduate School of Mathematical Sciences, The University of Tokyo, 3-8-1 Komaba, Meguro-ku, Tokyo 153-8914, Japan}
\email{suzukit@ms.u-tokyo.ac.jp}

\date{\today}

\begin{document}

\begin{abstract}
We prove that a knot $K \subset S^3$ is the unknot if and only if the $\Z/2$-equivariant Seiberg--Witten Floer cohomology of the double cover branched along the $(2, 1)$-cable of $K$ has no $\mathbb{F}_2[Q]$-torsion. The proof combines the Thom--Gysin exact sequence, the Montesinos trick, and the genus bound for $L$-space surgeries.
\end{abstract}

\maketitle

\section{Introduction}
A fundamental question about a knot invariant is whether it detects the unknot. 
For example, whether the Jones polynomial detects the unknot is a famous open problem; it is posed as Problem 1.24 in \cite{Baykur-Kirby-Ruberman2026K3}.
Heegaard/monopole/instanton knot Floer homology detect the Seifert genus 
\cite[Theorem~1.2]{OzsvathSzabo2004Genus}\cite[Lemmas~5.8 and~5.9]{KronheimerMrowka2010Excision}\cite[Proposition~7.16]{KronheimerMrowka2010Excision} and in particular detects the unknot. 
Using singular instanton Floer theory, Kronheimer--Mrowka \cite{Kronheimer-Mrowka2011Khovanov} proved Khovanov homology detects the unknot, and Dowlin \cite{Dowlin2024spectral} gave another proof using knot Floer homology.

From $\Z/2$-equivariant Seiberg--Witten Floer cohomology of the double branched cover of a knot $K \subset S^3$, the first author and Taniguchi \cite{Iida-Taniguchi2024monopoles} introduced an $\mathbb{F}_2[Q]$-module 
\[
\mathcal{H}_Q(K):=\widetilde{H}^{*}_{\mathbb{Z}/2}(SWF(\Sigma_2(K),\s_0); \mathbb{F}_2),
\]
where $\s_0$ is the $\Spinc$ structure induced by the unique Spin structure, and $\deg(Q)=1$. 
Here, the $\Z/2$-action on the Seiberg--Witten Floer spectrum is induced by the covering involution
on $\Sigma_2(K)$ together with its complex-linear lift to the spinor bundle (see \cite[Remark~3.1]{Iida-Taniguchi2024monopoles}).

Unfortunately, this invariant does not detect the unknot:
\begin{proposition}
Let $U$ be the unknot and $4_1$ the figure eight knot.
Then we have a graded isomorphism
\[
\mathcal{H}_Q(U)\cong \mathcal{H}_Q(4_1).
\]
\end{proposition}
\begin{proof}
Both $U$ and $4_1$ are quasi-alternating knots, and hence their double branched covers are $\mathbb{F}_2$-coefficient $L$-spaces.
Moreover, their signatures are
\[
\sigma(U)=\sigma(4_1)=0.
\]
It follows from \cite[Theorem~4.5]{Iida-Taniguchi2024monopoles}
that 
\[
\mathcal{H}_Q(U)\cong \mathcal{H}_Q(4_1)\cong \mathbb{F}_2[Q],
\] where $1\in \mathbb{F}_2[Q]$ has grading zero.
\end{proof}

However, the main result of this paper is that the knot invariant
\[
\mathcal{H}_Q(K_{2, 1})
\]
detects the unknot, where $K_{2, 1}$ is the $(2, 1)$-cable of the knot $K\subset S^3$.
Let us denote the torsion part of $\mathcal{H}_Q(K)$ by
\[
\mathcal{T}_Q(K):=\{x\in \mathcal{H}_Q(K): Q^n x=0 \text{ for some }n\geq 1\}.
\]
\begin{theorem}\label{thm:unknot-detect-conditions}
For every knot $K \subset S^3$, the following are equivalent:
\begin{enumerate}
\item $K$ is the unknot;
\item $\mathcal{T}_Q(K_{2, 1})=0$;
\item $\mathcal{H}_Q(K_{2, 1})\cong \mathbb{F}_2[Q]$ as ungraded $\mathbb{F}_2[Q]$-modules.
\end{enumerate}
\end{theorem}
One can also introduce a numerical invariant that detects the unknot.
Define
\[
N^{\cab}_Q(K)=\dim_{\mathbb{F}_2}(\mathcal{T}_Q(K_{2, 1})/Q\mathcal{T}_Q(K_{2, 1})).
\]
\begin{theorem}\label{thm:n-cab}
If $K$ is nontrivial, then $N^{cab}_Q(K)\geq 2$.
On the other hand, $N^{\cab}_Q(U)=0$.
Moreover, for any knot $K \subset S^3$, 
\[
N^{cab}_Q(K)=\frac{\dim_{\mathbb{F}_2}\widehat{HF}(S^3_{+1}(K\# K^r); \mathbb{F}_2)-1}{2}.
\]
Here, $K^r$ denotes $K$ with its  orientation reversed, not mirror image.
\end{theorem}

We assume that all manifolds are smooth, compact and oriented. 

\subsection*{Acknowledgements}
NI was supported by JSPS KAKENHI Grant Number JP26K16990 (Grant-in-Aid for Early-Career Scientists) and by the Start-up Fund of the Kavli Institute for the Physics and Mathematics of the Universe (Kavli IPMU), The University of Tokyo. 
TS was supported by JSPS KAKENHI, Grant Numbers JP26KJ0062 (Grant-in-Aid for JSPS Fellows) and JP26K16991 (Grant-in-Aid for Early-Career Scientists). 

The authors used ChatGPT to help clarify the proof of Lemma~\ref{lem:dimension}, conduct exploratory bibliographic
searches, and check English. The authors independently checked all references and
verified all mathematical arguments and computations, and take full responsibility for the content of
the manuscript.

\section{The Thom-Gysin sequence and a dimension formula}
Let $K \subset S^3$ be a knot.
By the rank one theorem~\cite[Theorem 3.8]{Iida-Taniguchi2024monopoles}, 
\[
\operatorname{rank}_{\mathbb{F}_2[Q]}\mathcal{H}_Q(K)=1.
\]
Therefore the structure theorem over $\mathbb{F}_2[Q]$ gives a decomposition
\[
\mathcal{H}_Q(K)\cong \mathbb{F}_2[Q]\langle q\rangle \bigoplus^N_{i=1} \mathbb{F}_2[Q]/(Q^{a_i})\langle d_i \rangle, \qquad a_i \geq 1,
\]
for some $q, d_i \in \Q$. 
Here $\mathbb{F}_2[Q]\langle d\rangle$ means that its generator has degree $d$.
Define
\[
N_Q(K):=N=\dim_{\mathbb{F}_2}(\mathcal{T}_Q(K)/Q\mathcal{T}_Q(K) )=\dim_{\mathbb{F}_2}\ker(Q: \mathcal{T}_Q(K)\to \mathcal{T}_Q(K)).
\]
\begin{lemma}\label{lem:dimension}
For any knot $K \subset S^3$, we have
\[
\dim_{\mathbb{F}_2}\widetilde{H}^*(SWF(\Sigma_2(K), \s_0); \mathbb{F}_2)= 1+2N_Q(K).
\]
In particular, 
\[
\mathcal{T}_Q(K)=0 \quad \Longleftrightarrow\quad \dim_{\mathbb{F}_2}\widetilde{H}^*(SWF(\Sigma_2(K), \s_0); \mathbb{F}_2)=1.
\]
\end{lemma}
\begin{proof}
By the exactness of the Thom--Gysin sequence used in the proof of Theorem 4.9 of \cite{Iida-Taniguchi2024monopoles}, we obtain 
\[
\cdots\to \mathcal{H}^{k-1}_Q(K) \xrightarrow{Q}\mathcal{H}^{k}_Q(K)
\to
\widetilde{H}^k(SWF(\Sigma_2(K), \s_0); \mathbb{F}_2)\to \mathcal{H}^{k}_Q(K) \xrightarrow{Q} \cdots. 
\]
By exactness, we obtain the short exact sequence
\[
0\to\mathcal{H}^{k}_Q(K)/ Q\mathcal{H}^{k-1}_Q(K)\to 
\widetilde{H}^k(SWF(\Sigma_2(K), \s_0); \mathbb{F}_2)\to \ker(Q:\mathcal{H}^{k}_Q(K)\to \mathcal{H}^{k+1}_Q(K))\to 0.
\]
Since
\[
\dim_{\mathbb{F}_2}(\mathcal{H}_Q(K)/ Q\mathcal{H}_Q(K))=1+N_Q(K)
\]
and
\[
\dim_{\mathbb{F}_2}\ker(Q:\mathcal{H}_Q(K)\to \mathcal{H}_Q(K))=N_Q(K),
\]
we have
\[
\dim_{\mathbb{F}_2}\widetilde{H}^*(SWF(\Sigma_2(K), \s_0); \mathbb{F}_2)=1+2N_Q(K).
\]
\end{proof}

A rational homology $3$-sphere $Y$ is called an {\it $L$-space} if 
\[
\dim_{\mathbb{F}_2}\widehat{HF}(Y;\mathbb{F}_2)
=|H_1(Y;\mathbb{Z})|.
\]
Throughout this paper, we use the term $L$-space with this coefficient convention. 

\begin{corollary}\label{cor:tql-iff-l-space}
Let $K \subset S^3$ be a knot.
If $\Sigma_2(K)$ is an integral homology $3$-sphere, then $\mathcal{T}_Q(K)=0$ if and only if $\Sigma_2(K)$ is an $L$-space. 
\end{corollary}

\begin{proof}
Since $\Sigma_2(K)$ is an integral homology $3$-sphere, it has a unique $\mathrm{Spin}^c$ structure. 
Hence,
\[
|H_1(\Sigma_2(K);\mathbb{Z})|=1.
\]
The desired equivalence follows from Lemma~\ref{lem:dimension}.
\end{proof}

\section{The Montesinos trick and Floer homology computations}
In this section, we prove lemmas needed later. 
We write $X_1 \cong X_2$ to mean that the manifolds $X_1$ and $X_2$ are diffeomorphic, and $K_1 \sim K_2$ to mean that the knots $K_1$ and $K_2$ in $S^3$ are isotopic. 
For a knot $K \subset S^3$, let $Y_K := \Sigma_2(K_{2,1})$.

\begin{lemma}\label{lem:Montesinos-trick-homeo}
For any knot $K \subset S^3$, we have
\[
Y_K \cong S^3_{+1}(K\#K^r).
\]
In particular, $Y_K$ is an integral homology $3$-sphere.
\end{lemma}

\begin{proof}
The diffeomorphism $Y_K \cong S^3_{+1}(K\#K^r)$ follows from the Montesinos trick. 
The details, including the involution on the branched cover, are described in \cite[Section~2]{Dai-Sungkyung-Mallick-Park-Stoffregen2024cable}. 

The second assertion follows from the facts that $H_1(S^3_{+1}(J))=0$ for any knot $J \subset S^3$ and that any closed $3$-manifold $Y$ with $H_1(Y)=0$ is an integral homology $3$-sphere. 
\end{proof}

For any knot $K \subset S^3$, we have
\begin{equation}\label{eq:genus-connected-sum}
g(K\#K^r)
=g(K)+g(K^r)
=2g(K),
\end{equation}
where $g(J)$ denotes the Seifert genus of a knot $J$.

For a knot $J \subset S^3$ and a positive integer $n$, if the $3$-manifold $S^3_n(J)$ obtained by $n$-surgery along $J$ is an $L$-space, then
\begin{equation}\label{eq:L-space-genus-bound}
n\geq 2g(J)-1
\end{equation}
(see \cite[Corollary 8.5]{Kronheimer-Mrowka-Ozsvath-Szabo2007monopoles} and \cite[Section 1.1]{Greene2015L-space}).

\begin{lemma}\label{lem:L-iff-trivial}
The manifold $Y_K$ is an $L$-space if and only if $K$ is the unknot.
\end{lemma}

\begin{proof}
Suppose that $Y_K$ is an $L$-space. By Lemma~\ref{lem:Montesinos-trick-homeo}, we have $Y_K\cong S^3_{+1}(K\#K^r)$. 
It follows from \eqref{eq:genus-connected-sum} and \eqref{eq:L-space-genus-bound} that
\[
4g(K)-1 = 2g(K\#K^r)-1 \leq 1.
\]
Hence, $g(K)\leq 1/2$. 
Since the Seifert genus $g(K)$ is a nonnegative integer, we obtain $g(K)=0$. Therefore, $K$ is the unknot.

Conversely, suppose that $K$ is the unknot $U$. 
Then $U_{2,1} \sim U$, and hence $Y_U \cong S^3$. 
Since $S^3$ is an $L$-space, so is $Y_U$.
\end{proof}

\begin{lemma}\label{lem:dim-geq-5}
For any nontrivial knot $K \subset S^3$, we have
\[
\dim_{\mathbb{F}_2}
\widehat{HF}(S^3_{+1}(K\#K^r);\mathbb{F}_2)
\geq 5.
\]
\end{lemma}

\begin{proof}
Set $J=K\#K^r$. 
By the exactness of the surgery exact triangle for $\widehat{HF}$ \cite[Theorem 9.16]{Ozsvath-Szabo2004holomorphic}, 
\[
\cdots \to
\widehat{HF}(S^3;\mathbb{F}_2) \to
\widehat{HF}(S^3_0(J);\mathbb{F}_2) \to
\widehat{HF}(S^3_{+1}(J);\mathbb{F}_2) \to
\widehat{HF}(S^3;\mathbb{F}_2) \to
\cdots,
\]
we obtain
\[
\dim_{\mathbb{F}_2}
\widehat{HF}(S^3_0(J);\mathbb{F}_2)
\leq
\dim_{\mathbb{F}_2}\widehat{HF}(S^3;\mathbb{F}_2)
+
\dim_{\mathbb{F}_2}\widehat{HF}(S^3_{+1}(J);\mathbb{F}_2).
\]
Since $K$ is nontrivial, we have $2g(K)>1$. 
It follows from the proof of Theorem 1.1 in \cite{Hedden2009Khovanov} that 
\[
\dim_{\mathbb{F}_2}
\widehat{HF}(S^3_0(J);\mathbb{F}_2)
\geq 6. 
\]
Combining these inequalities with $\dim_{\mathbb{F}_2}\widehat{HF}(S^3;\mathbb{F}_2)=1$, we obtain 
\[
\dim_{\mathbb{F}_2}\widehat{HF}(S^3_{+1}(J);\mathbb{F}_2)
\geq
\dim_{\mathbb{F}_2}\widehat{HF}(S^3_0(J);\mathbb{F}_2)
-
\dim_{\mathbb{F}_2}\widehat{HF}(S^3;\mathbb{F}_2)
\geq 5.
\]
Since $J=K\#K^r$, the desired inequality follows. 
\end{proof}

\begin{lemma}\label{lem:HF=HM}
For any closed oriented $\Spinc$ rational homology $3$-sphere $(Y, \s)$, we have
\[
\widehat{HF}(Y, \s;\mathbb{F}_2)
\cong
\widetilde{H}^*(SWF(Y,\s); \mathbb{F}_2).
\]
\end{lemma}

\begin{proof}
The isomorphism follows from the equivalence $HF=HM$ established in \cite{Kutluhan-Lee-Taubes2020I}, together with the comparison between monopole Floer
homology and the Seiberg--Witten Floer spectrum \cite{Lidman-Manolescu2017equivalence}.
\end{proof}

\section{Proofs of the main theorems}
In this section we prove Theorems \ref{thm:unknot-detect-conditions} and \ref{thm:n-cab}.

\begin{proof}[Proof of Theorem~\ref{thm:unknot-detect-conditions}]
\noindent
\textup{(1) $\implies$ (3):}
Suppose that $K$ is the unknot $U$. Then $K_{2,1} \sim T(2,1) \sim U$. 
Hence, we have $\mathcal{H}_Q(K_{2,1})\cong \mathcal{H}_Q(U)\cong \mathbb{F}_2[Q]$ as ungraded $\mathbb{F}_2[Q]$-modules. 

\medskip
\noindent
\textup{(3) $\implies$ (2):}
Suppose that $\mathcal{H}_Q(K_{2,1})\cong \mathbb{F}_2[Q]$ as ungraded $\mathbb{F}_2[Q]$-modules. 
Since $\mathcal{H}_Q(K_{2,1})$ has no $Q$-torsion, we obtain $\mathcal{T}_Q(K_{2,1})=0$.

\medskip
\noindent
\textup{(2) $\implies$ (1):}
Suppose that $\mathcal{T}_Q(K_{2,1})=0$.
By Corollary~\ref{cor:tql-iff-l-space}, $Y_K=\Sigma_2(K_{2,1})$ is an $L$-space. 
It follows from Lemma~\ref{lem:L-iff-trivial} that $K$ is the unknot. 

Thus, the three conditions (1), (2), and (3) are equivalent.
\end{proof}

\begin{proof}[Proof of Theorem~\ref{thm:n-cab}]
By Lemma~\ref{lem:Montesinos-trick-homeo}, we have $Y_K=\Sigma_2(K_{2,1}) \cong S^3_{+1}(K\#K^r)$. 
Since $Y_K$ is an integral homology $3$-sphere, it has a unique $\mathrm{Spin}^c$ structure $\mathfrak{s}_0$.

By Lemmas~\ref{lem:dimension} and \ref{lem:HF=HM}, we have
\[
\dim_{\mathbb{F}_2}\widehat{HF}(\Sigma_2(K_{2,1}); \mathbb{F}_2)
=
\dim_{\mathbb{F}_2}\widetilde{H}^*(SWF(\Sigma_2(K_{2,1})); \mathbb{F}_2)
=
1+2N_Q^{\mathrm{cab}}(K).
\]
Therefore,
\[
N_Q^{\mathrm{cab}}(K)
=
\dfrac{\dim_{\mathbb{F}_2}\widehat{HF}(S^3_{+1}(K\#K^r); \mathbb{F}_2)-1}{2}.
\]
Moreover, since $\dim_{\mathbb{F}_2}\widehat{HF}(S^3_{+1}(U\#U^r); \mathbb{F}_2)
=
\dim_{\mathbb{F}_2}\widehat{HF}(S^3; \mathbb{F}_2)
=1$, we have $N_Q^{\mathrm{cab}}(U)=0$. 
On the other hand, for any nontrivial knot $K$, Lemma~\ref{lem:dim-geq-5} gives 
\[
\dim_{\mathbb{F}_2}\widehat{HF}(S^3_{+1}(K\#K^r); \mathbb{F}_2)\geq 5. 
\]
Hence, $N_Q^{\mathrm{cab}}(K)\geq 2$.
\end{proof}
\bibliographystyle{amsalpha}
\bibliography{bibliography}

\end{document}